\documentclass[12pt]{article}
\usepackage[pdftex]{graphicx}
\usepackage{amsfonts}
\usepackage{amsmath}

\newtheorem{theorem}{Theorem}[section]
\newtheorem{lemma}[theorem]{Lemma}
\newtheorem{corollary}[theorem]{Corollary}

\newtheorem{example}[theorem]{Example}

\newtheorem{remark}[theorem]{Remark}

\numberwithin{equation}{section}

\newcommand*{\os}{\widetilde}
\newcommand*{\pd}[2]{\frac{\partial #1}{\partial #2}}

\begin{document}

\title{On non-local reflection for  elliptic equations of the second
order in $\mathbb{R}^2$\\
(the Dirichlet condition)}
\author{Tatiana Savina \\
Department of Mathematics\\
321 Morton Hall, Ohio University\\
Athens, OH 45701, U.S.A.}

\maketitle


\date{}

\begin{abstract}
Point-to-point reflection holding for harmonic functions subject to the Dirichlet or Neumann conditions on an analytic curve in the
plane almost always fails for solutions to more general elliptic
equations. We develop a non-local, point-to-compact set, formula for
reflecting a solution of an analytic elliptic partial differential
equation across a real-analytic curve on which it satisfies the
Dirichlet conditions. We also discuss the special cases when the formula reduces to the point-to-point forms.

\end{abstract}




\section{Introduction}

Schwarz symmetry principle is one of the celebrated tools in analysis and mathematical physics that has been attracting attention of many mathematicians \cite{as}-- \cite{john}, \cite{laguna}--\cite{lopez}, \cite{poritsky}--\cite{shapiro}.  From the point of view of applications it is important to have an explicit reflection formula for a specific problem (\cite{duffin}, \cite{farwig}, \cite{poritsky}).
One of the open questions is the following:  for what partial differential equations, boundary conditions and spatial dimensions such a formula exists and what is the structure of this formula, in other words, whether it is a point to point formula (see, for example \cite{ek}) or it
 has a more complicated structure, for example,
 a point to a finite set \cite{lopez} or a point to a continuous set  (see, for example, \cite{bs} and references therein).  

In this paper, we derive a  reflection formula for solutions of 
elliptic equations in $\mathbb{R}^2$ with respect to  a non-singular
real-analytic curve and study the obtained formula. We  call this
formula  non-local, since unlike the classical point to point reflection (see the Theorem \ref{Th1} below)  
this is a point to compact set reflection, 
generalizing  the following celebrated Schwarz reflection
principle for harmonic functions.
\begin{theorem}\label{Th1} (\cite{laguna} Chapter 9, p.~51; \cite{shapiro} Chapter 1, p.~4).
Let $\Gamma =\{(x, y): f(x,y)=0\} \subset \mathbb{R}^2$ be a non-singular real-analytic curve
and $P'\in \Gamma$. Then, there exists a neighborhood $U$ of $P'$ and an
anti-conformal mapping $R:U\to U$ which is identity on $\Gamma$,
permutes the components $U_1, U_2$ of $U\setminus \Gamma$ and relative to
which any harmonic function $u(x, y)$ defined near $\Gamma $ and
vanishing on $\Gamma $ is odd; i.e.,
\begin{equation} \label{E:0.1}
         u(x_0 ,y_0 ) =-u(R(x_0 ,y_0 ))
\end{equation}
for any point $(x_0, y_0)$ sufficiently close to $\Gamma$.
Note that if the point $(x_0 ,y_0)\in U_1$, then the
``reflected'' point $R(x_0 ,y_0)\in U_2$.
\end{theorem}

Here the mapping $R$ can be  described  by considering a complex domain $U_{\mathbb{C}}$ in the
space $\mathbb{C}^2$,  such that $U_{\mathbb{C}}\cap \mathbb{R}^2=U$, to which the function $f$, defining the curve
$\Gamma$, is continued analytically.
After the
transformation of the variables, $z=x+iy$, $\zeta =x-iy,$ the equation of the
complexified curve $\Gamma _{\mathbb{C}}$ can be rewritten in the
form
\begin{equation}\label{E:999}
f\left ( \frac{z+\zeta }{2}, \, \frac{z-\zeta }{2i} \right ) =0.
\end{equation}
If $grad \, f(x,y)\ne 0$ on $\Gamma$, \eqref{E:999} can be solved
with respect to $z$ or $\zeta$; the corresponding solutions we denote
as $\zeta =S(z)$ and $z=\overset\sim{S}(\zeta )$. The function $S(z)$ is
called the {\it Schwarz function} of the curve $\Gamma$
\cite{davis}, Chapter 5, p.~21. The mapping $R$ is given by
\begin{equation} \label{E:2}
R(x ,y )=R(z)=\overline{S(z )}.
\end{equation}

Formula \eqref{E:0.1} has been generalized to cover several other
situations including the Helmholtz equation and wave equation, and the polyharmonic functions (see, for example,  \cite{s99}, \cite{as}, \cite{lopez} and references therein). The purpose of this paper is to obtain an explicit
reflection formula for solutions to the elliptic equation 
\begin{equation}\label{intro1}
{ L}u\equiv \Delta _{x,y}u+a(x,y)\pd{u}{x}+b(x,y)\pd{u}{y}+c(x,y)u=0
\end{equation}
with respect to  a real analytic curve in $\mathbb{R}^2$, where the solution vanishes,
 and to investigate the properties
of the mapping induced by this formula.
Here $a(x,y)$, $b(x,y)$ and $c(x,y)$ are real-analytic functions in
the  domain  $U\subset\mathbb{R}^2$.

 In what follows,  
a formula, expressing the value of a function $u(x,y)$ at an
arbitrary point $(x_0 ,y_0)\in U_1$ in terms of its values at points
in $U_2$, is called a reflection formula. It is more often an integro-differential operator than a  point-to-point reflection \eqref{E:0.1}, which seems
to be quite rare for solutions of partial differential equations. In
particular, for solutions of the Helmholtz equation $(\Delta _{x,y}
+\lambda ^2 )u(x,y)=0$ vanishing on a curve $\Gamma$, point-to-point
reflection holds
only when $\Gamma$ is a line, while for harmonic functions in
$\mathbb{R}^3$ it holds only when $\Gamma$ is either a plane or a
sphere  \cite{ek}, \cite{khavinsons}. The paper by P.~Ebenfelt and
D.~Khavinson \cite{ek} is devoted to further study of point-to-point
reflection for harmonic functions. There it was shown that
point-to-point reflection in the sense of the Schwarz reflection
principle for $n>2$ is very rare in $\mathbb{R}^n$ when $n$ is even, and
that it never holds when $n$ is odd, unless $\Gamma$ is a
sphere or a hyperplane. Reflection properties of solutions of the
Helmholtz equation have also been considered in \cite{ebenfelt},
\cite{sss}, \cite{s99}. Two later papers are devoted to derivation
of non-local formulas for Helmholtz equation subject to Dirichlet and
Neumann conditions respectively. 
Recently a  reflection formula for harmonic functions subject to the Robin condition, $\alpha\, \partial _{n} u
+\beta\, u=0$, on a real-analytic curve was derived in \cite{bs}, and it was
shown  that the obtained (non-local) formula reduces to well-known point-to-point
reflection laws corresponding to the Dirichlet and  Neumann boundary
conditions when one of the coefficients, $\alpha$ or $\beta$,
vanishes.

The structure of the paper is as follows: in  Section \ref{sec2} we describe some preliminaries; in  Section \ref{sec3} we formulate the main theorem, which is proven in  Section \ref{sec4}. Conclusions and the special cases, when the point to point reflections hold, are discussed in  Section \ref{sec5}.

\section{Preliminaries}\label{sec2}

We are starting this section by recalling a classical B. Riemann's result for hyperbolic equations (see \cite{hadamard} Chapter 2, p.~65 or \cite{Gar} Chapter 4, p.~127 for detailed explanations, here we follow a short version \cite{laguna} Chapter 9, p.~55).

Consider a hyperbolic differential equation with entire coefficients
\begin{equation}\label{hyperb}
 { H}u \equiv \frac{\partial ^2 u}{\partial
x\partial y} +a(x,y)\pd{u}{x} +b(x,y)\pd{u}{y} +c(x,y)u=0,
\end{equation}
its adjoint equation is 
\begin{equation}\label{hyperba}
 { H}^* u\equiv \frac{\partial ^2 u}{\partial
x\partial y} -\pd{(au)}{x} -\pd{(bu)}{y} +cu=0.
\end{equation}
The Riemann function ${\mathfrak R}_H (x,y;x_0,y_0)$ of operator $ H$ is defined as the solution to the Goursat problem:
\begin{equation}\label{riemann}
\begin{cases}
& { H}^* {\mathfrak R}_H =0 \quad \text{ near } (x_0,y_0),\\
& {\mathfrak R}_H (x_0,y;x_0,y_0)= \exp \Bigl \{\int\limits _{y _0} ^{
y} a(x _0,\tau )d\tau \Bigr \}, \\
& {\mathfrak R}_ H(x,y_0;x_0,y_0)= \exp \Bigl \{\int\limits _{x _0} ^{
x} b(t,y _0 )dt \Bigr \}.
\end{cases}
\end{equation}
Note that ${\mathfrak R}_H$ is an entire function of all four variables, moreover ${\mathfrak R}_H (x,y;x_0,y_0)={\mathfrak R}_{H^*} (x_0,y_0;x,y)$,  ${\mathfrak R}_H (x_0,y_0;x_0,y_0)=1$ and
the following Riemann's lemma holds.

\begin{lemma} Let $\Gamma := \{ (x,y)| \, y=s(x)\}$ be a non-characteristic with respect to $H$ real-analytic curve that divides a domain $U\subset \mathbb{R}^2$ in two connected components $U_1$ and $U_2$, and $u(x,y)$ be a solution of (\ref{hyperb}) near $\Gamma$. For all points $P(x_0,y_0)\in U$ sufficiently close to $\Gamma$
we have
\begin{equation}\label{hypref}
u(P)=\frac{1}{2}u(M){\mathfrak R}_H(M)+\frac{1}{2}u(N){\mathfrak R}_H(N)-\int\limits _M ^N ({\mathfrak U} dy-{\mathfrak V} dx),
\end{equation}
where $M=(s^{-1}(y_0),y_0)$, $N=(x_0,s(x_0))$ and
\begin{eqnarray}\nonumber
{\mathfrak U}=a{\mathfrak R}_H+\frac{1}{2}{\mathfrak R}_H\pd{u}{y}-\frac{1}{2}\pd{{\mathfrak R}_H}{y}u,\\
{\mathfrak V}=b{\mathfrak R}_H+\frac{1}{2}{\mathfrak R}_H\pd{u}{x}-\frac{1}{2}\pd{{\mathfrak R}_H}{x}u.\nonumber
\end{eqnarray}
\end{lemma}

If in addition solution to the equation $Hu=0$  vanishes on $\Gamma$, formula (\ref{hypref}) reduces to
\begin{equation}\label{hypref1}
u(P)=\frac{1}{2}\int\limits _M ^N {\mathfrak R}_H(\pd{u}{x}dx-\pd{u}{y} dy).
\end{equation}

\begin{remark}
For the wave equation, $a=b=c=0$ in (\ref{hyperb}), the Riemann function equals 1 identically.
Consider a point $P(x_0,y_0)\in U_1$ and
  a solution of the wave equation   vanishing on $\Gamma$. Let's
 allow the path of integration $MN$ in (\ref{hypref1}) to degenerate to a pair of segments (in $U_2$) of  a vertical and horizontal characteristics through points
$M$ and  $N$, which intersect at a point $Q(s^{-1}(y_0),s(x_0))\in U_2$. Then  formula (\ref{hypref1}) becomes
$$
u(P)=-u(Q).
$$

Since the points $P$ and $Q$ are located on the opposite sides of the curve $\Gamma$, the later formula  states a  point to point reflection law for the wave equation.
\end{remark}
\begin{remark}
If for a solution of the wave equation vanishing on $\Gamma$ we allow the path of integration $MN$ in (\ref{hypref1}) to degenerate to a polygonal line  consisting of vertical and horizontal segments 
with vertices $M=Q_1, Q_2, Q_3,\dots , Q_n=N$ such that $Q_{2k+1}\in\Gamma$ and $Q_{2k}\in U_2\setminus\Gamma$,
$k=1,\dots ,(n-1)/2$,
then a version of a point to finite set reflection will be obtained \cite{lopez}
$$
u(P)=-\sum\limits _{k=1}^{(n-1)/2} u(Q_{2k}) 
$$
(see \cite{lopez} for other examples of point to finite set formulas).
\end{remark}

If we consider the elliptic equation (\ref{intro1}) in the complex domain $U_{\mathbb{C}}\subset\mathbb{C}^2$, then the equation and its adjoint in characteristic variables $(z,\zeta )$ become similar to the hyperbolic equation (\ref{hyperb}) and its adjoint (\ref{hyperba})
\begin{equation}\label{eqc}
 { L} _{\mathbb C}{u}\equiv \frac{\partial ^2 u}{\partial
z\partial \zeta} +A\pd{u}{z} +B\pd{u}{\zeta} +Cu=0,
\end{equation}
\begin{equation}\label{eqca}
{ L}^* _{\mathbb C}{u}\equiv \frac{\partial ^2 u}{\partial
z\partial \zeta} -\pd{(Au)}{z}-\pd{(Bu)}{\zeta} +Cu=0,
\end{equation}
where the coefficients in (\ref{intro1}) are replaced with
$$
A(z,\zeta )=\frac{1}{4}\bigl [a(x,y)+ib(x,y) \bigr ],
\quad
B(z,\zeta )=\frac{1}{4}\bigl [a(x,y)-ib(x,y) \bigr ],
$$
$$
C(z,\zeta )=\frac{1}{4}c(x,y).
$$
Analogously, the Riemann function of $L$ is defined as the solution to the Goursat problem in $\mathbb{C}^2$:
\begin{equation}\label{riemann}
\begin{cases}
& { L}^* _{\mathbb C}{\mathfrak R}\equiv
\frac{\partial ^2 }{\partial z\partial \zeta}{\mathfrak R}
-\pd{}{z}(A{\mathfrak R})
-\pd{}{\zeta}(B{\mathfrak R}) +C{\mathfrak R}=0,\\
& {\mathfrak R}_ {|_{z=z_0}}= \exp \Bigl \{\int\limits _{\zeta _0} ^{
\zeta} A(z _0,\tau )d\tau \Bigr \}, \\
& {\mathfrak R}_ {|_{\zeta =\zeta _0}}= \exp \Bigl \{\int\limits _{z _0} ^{
z} B(t,\zeta _0 )dt \Bigr \}.
\end{cases}
\end{equation}

By a  {\it  fundamental
solution} of operator $ L$ we understand a solution of the equation ${
L}^*G(x_0,y_0,x,y)=\delta (x_0,y_0)$, where ${ L}^*$ is the adjoint to $ L$
differential operator. Thus, function $G$ written in the characteristic variables $z=x+iy$ and
$\zeta =x-iy$ is a solution to the  equation
\begin{equation}\label{ur}
 L^*_{\mathbb C}u=\frac{\partial ^2u}{\partial z\partial \zeta}
-\frac{\partial Au}{\partial z}-\frac{\partial Bu}{\partial
\zeta}+Cu=\delta (z_0,\zeta _0).
\end{equation}
The following formula (see \cite{hadamard} Chapter 3, p.~72) shows that the Riemann function is a factor of the logarithm in an expression for the fundamental solution of the operator $L_{\mathbb C}$
\begin{equation}\label{r1}
G(z,\zeta ; z_0,\zeta _0)=-\frac{1}{4\pi}{ \mathfrak R}(z, \zeta ; z_0 ,\zeta
_0 )\ln [(z-z_0 )(\zeta -\zeta _0)]+g_0(z,\zeta ,
z_0,\zeta _0),
\end{equation}
where $g_0(z,\zeta , z_0,\zeta _0)$ is an entire function.

Note, that the
fundamental solution exists (see \cite{john1} Chapter 3, p.~50) and is uniquely determined   up to the kernel of operator $
L$. 

There are different representations of the fundamental solution, for example,  \cite{colton}; \cite{john1} Chapter 3, p.~76;  \cite{john2}. However, for what follows we need  
a special representation  as a sum of two functions, each of those has the logarithmic singularity on a single characteristic in $\mathbb{C}^2$. This representation  is given by the following theorem.
\begin{theorem}\label{L1} \cite{fundsol} There exist a fundamental solution of $ L$ that
can be represented in
the form
\begin{equation}\label{G}
G=-\frac{1}{4\pi}(G_1+G_2),
\end{equation}
\begin{equation}\label{Gj}
G_j=\sum\limits _{k=0}^{\infty}\alpha _k^j(x_0,y_0;x,y)f_k(\psi
_j),\qquad j=1,2\, ,
\end{equation}
\begin{equation}\label{fk}
f_k(\xi )=
\begin{cases}
&(-1)^{-k-1}(-k-1)!\xi ^k,\qquad k\le -1,\\
&\frac{\xi ^k}{k!}(\ln \xi -C_k),\qquad k=0,1,...\, ,
\end{cases}
\end{equation}
\begin{equation}\notag
 C_0=0,\qquad C_k=\sum\limits _{l=1}^k \frac{1}{l},\qquad
k=1,2,...\, ,
\end{equation}
\begin{equation}
\psi _1=(x-x_0)+i(y-y_0 )=z-z_0, \qquad \psi _2=(x-x_0)-i(y-y_0
)=\zeta -\zeta _0.
\end{equation}
Here the coefficients  $\alpha _k^j$ are uniquely determined by 
  recursive transport equations
\begin{equation}\label{transp}
\begin{aligned}
& {\mathfrak L}\alpha _0^j=0,\qquad {\mathfrak L}\alpha
_{k+1}^j=- L^*_{\mathbb
C}\alpha _k^j,  \\
&{\mathfrak L}=\pd{\psi _j}{z}\cdot \Bigl [ \pd{}{\zeta}-A\Bigr ]
+\pd{\psi _j}{\zeta}\cdot \Bigl [\pd{}{z}-B\Bigr ] 
\end{aligned}
\end{equation}
subject to the  initial conditions
\begin{equation}\label{initcond}
\begin{cases}
& {\alpha _0^1}_ {|_{\zeta =\zeta _0}}= \exp \Bigl \{\int\limits _{z
_0} ^{ z} B(t,\zeta _0 )dt \Bigr \},\,\, {\alpha _0^1}_ {|_{\zeta
=\zeta _0}}=0, \,\, k=1,2,\dots \, , \\
 & {\alpha _0^2} _ {|_{z=z_0}}=
\exp \Bigl \{\int\limits _{\zeta _0} ^{ \zeta} A(z _0,\tau )d\tau
\Bigr \},\quad {\alpha _k^2} _ {|_{z=z_0}}=0, \quad k=1,2,\dots \, .
\end{cases}
\end{equation}
\end{theorem}
Note that  \eqref{transp} and \eqref{initcond}, in particular, imply 
\begin{equation}
\alpha _0^1=\exp \Bigl (\int\limits _{\zeta _0}^{\zeta}A(z,\tau
)d\tau +\int\limits _{z _0}^{z}B(t,\zeta _0 )dt\Bigr ),\quad \alpha
_0^2=\exp \Bigl (\int\limits _{\zeta _0}^{\zeta}A(z_0,\tau )d\tau
+\int\limits _{z _0}^{z}B(t,\zeta  )dt\Bigr ).
\end{equation}
Taking into account (\ref{r1}), one can  interpret  $\alpha _k^j$ as coefficients in the following
series representations for the Riemann function \eqref{riemann} \cite{fundsol}:
\begin{equation}\label{anal}
{\mathfrak R}(z_0 ,\zeta _0 , z, \zeta )=\sum\limits _{k=0}^{\infty}
\alpha _k^1 (z_0,\zeta _0, z, \zeta ) \frac{(z-z_0)^k}{k!} =
\sum\limits _{k=0}^{\infty} \alpha _k^2 (z_0,\zeta _0, z, \zeta
)\frac{(\zeta-\zeta _0)^k}{k!}.
\end{equation}

\begin{remark}
For the Laplace equation,  $a=b=c=0$, and, therefore, $A=B=C=0$. Thus, $\alpha _0^1=\alpha _0^2=1$ and $\alpha _j^1=\alpha _j^2=0$, $j\ge 1$, and 
\begin{equation}
G_1^L  =\ln (z-z_0 ) \quad \text{and}\quad G_2^L
=\ln (\zeta-\zeta _0 )
\end{equation}
respectively, which leads to a standard fundamental solution:
\begin{equation}
G=-\frac{1}{4\pi}\ln  [ (z-z_0 )(\zeta-\zeta _0 )]= -\frac{1}{4\pi}\ln [(x-x_0)^2+(y-y_0)^2)].
\end{equation}
\end{remark}

\begin{remark}
In the case of the Helmholtz equation $a=b=0$ and $c=\lambda ^2$, that is, $\frac{\partial ^2 }{\partial
z\partial \zeta} u^H +\frac{\lambda ^2}{4}u^H =0$. Here $\lambda$
is a real number, functions $G_1$ and $G_2$ reduce to the form used in
\cite{sss}, \cite{s99}
\begin{equation}\label{fundH}
\begin{aligned}
G_1^H & =\sum\limits _{k=0}^{\infty} \frac{[-\lambda
^2(z-z_0)(\zeta -\zeta _0)]^k}{4^k (k!)^2}\Bigl ( \ln (z-z_0 ) -C_k\Bigr ), \\
G_2^H & =\sum\limits _{k=0}^{\infty} \frac{[-\lambda
^2(z-z_0)(\zeta -\zeta _0)]^k}{4^k (k!)^2}\Bigl ( \ln (\zeta -\zeta
_0 ) -C_k\Bigr )
\end{aligned}
\end{equation}
Summing up $G_1^H$ and $G_2^H$ and multiplying by $-\frac{1}{4\pi}$ one  obtains well-known fundamental
solution of the Helmholtz equation:
\begin{equation}
 -\frac{1}{4\pi}(G_1^H+G_2^H)=\frac{c+\ln \lambda /2}{2\pi}J_0 \Bigl ( \lambda \sqrt
{(z-z_0)(\zeta -\zeta _0)}\Bigr )-\frac{1}{4}N_0\Bigl ( \lambda
\sqrt {(z-z_0)(\zeta -\zeta _0)}\Bigr ),
\end{equation}
where $c$ is the Euler constant, and $J_0$ and $N_0$ are the Bessel
and the Neumann functions of zero order respectively.
\end{remark}

\section{ The main result  }\label{sec3}

Consider a solution of homogeneous linear elliptic differential
equation, written in its canonical form \cite{Gar} Chapter 5, p.~136 (with the
Laplace operator, $\Delta _{x,y}$, in the principal part), in a domain
$U\subset {\mathbb R}^2$  vanishing on an algebraic curve $\Gamma$,
\begin{equation}\label{pde}
  \begin{cases}
     & L u\equiv \Delta  _{x,y}u+a\pd{u}{x}+b\pd{u}{y}+cu=0\text{ near
}\Gamma ,\\
 &u(x,y)_{\mid _\Gamma }=0 ; a,\,b,\, c\,\,\text{ are real-analytic functions of
 }x,\,y.
\end{cases}
\end{equation}

\begin{theorem} \label{theorem}
Under the above assumptions, the following reflection formula holds in $U$:

\begin{equation}\label{formula}
 \begin{aligned}
\, & u(P) \,   =\,  -\,  c_0(P,\Gamma )\, u(Q)\,  +\\
& \frac{1}{2i}\int\limits   _{\Gamma }^Q \Bigl (
\bigl
\{
u\pd{V}{x}-V\pd{u}{x}-auV \bigr \} dy
      - \bigl \{u\pd{V}{y}-V\pd{u}{y}-buV \bigr \} dx
\Bigr ),
\end{aligned}
\end{equation}
where $P=(x_0 ,y_0 )$ and $Q=R(P)$ (see (\ref{E:2})), and the integral is computed along any curve joining $\Gamma$ with $Q$. Here
\begin{equation}\label{c0}
 \begin{aligned}
c_0(P,\Gamma ) &=\frac{1}{2}\Bigl \{ \exp \Bigl [\int\limits _{z_0} ^{{\os
S}(\zeta _0)} B(t,S(z_0))dt+\int\limits _{\zeta _0} ^{
S(z _0)} A(z_0 ,\tau)d\tau \Bigr ] \\
& +  \exp \Bigl [\int\limits _{\zeta _0} ^{ S(z _0)} A({\os S}(\zeta
_0),\tau )d\tau+\int\limits _{z _0} ^{ {\os S}(\zeta  _0)} B(t
,\zeta _0)dt \Bigr ] \Bigr \} ,
\end{aligned}
\end{equation}
where
$$
A(z,\zeta )=\frac{1}{4}\bigl [a(x,y)+ib(x,y) \bigr ],
\quad
B(z,\zeta )=\frac{1}{4}\bigl [a(x,y)-ib(x,y) \bigr ],
$$
$$
C(z,\zeta )=\frac{1}{4}c(x,y),\quad  V=V(x_0 ,y_0 ,x,y)=V_1 (x_0 ,y_0 ,x,y)-V_2 (x_0 ,y_0 ,x,y).
$$
Functions $V_j$ are solutions of the Cauchy-Goursat problems:
\begin{equation}\notag
\begin{cases}
\begin{aligned}
& { L}^*_{\mathbb{C}} V_j=0, && j=1,2\, ,\\
& {V_j}_ {|_{\Gamma _{\mathbb{C}}}}={\mathfrak R}_{|_{\Gamma _{\mathbb{C}}}},  &&  j=1,2\, , \\
& V_1= \exp \Bigl \{\int\limits _{\zeta _0} ^{
\zeta} A(\os{S}(\zeta ),\tau )d\tau +
\int\limits _{z _0} ^{z} B(t,\zeta  )dt\Bigr \},
&& \text{ on the char. } \os{l}_1=\{\os{S}(\zeta )=z_0\}\\
& V_2=     \exp \Bigl \{\int\limits _{\zeta _0} ^{
\zeta} A(z ,\tau )d\tau +
\int\limits _{z _0} ^{z} B(t,S(z)  )dt\Bigr \},
&& \text{ on the char. }\os{l}_2=\{{S}(z)=\zeta _0\},
\end{aligned}
\end{cases}
\end{equation}
where  ${ L}^*_{\mathbb{C}}$  is the adjoint operator  to $L_{\mathbb{C}}$ and
${\mathfrak R}(z_0 ,\zeta _0 , z,  \zeta )$ is the Riemann function of
$L$.

\end{theorem}

\section{Proof of the Theorem \ref{theorem}}\label{sec4}
\subsection{Sketch of the proof}\setcounter{equation}{0}
 We
begin with  Green's formula expressing a solution of  the  equation $Lu=0$
at a point $P$ via its values on a
contour $\gamma\subset U_1$ surrounding the point $P$  \cite{G}:
\begin{equation}\label{E:1.20}
    u(P)
=\int\limits_{\gamma} \omega [ u,G  ],
\end{equation}
where
\begin{equation}\label{E:1.20b}
   \omega  [  u,G ]=  \bigl \{ u\pd{G}{x}-G\pd{u}{x}-auG
\bigr \}\, dy
      - \bigl \{u\pd{G}{y}-G\pd{u}{y}-buG \bigr \}\, dx.
\end{equation}
Here  $G=G(x,y,x_0 ,y_0 )$ is an arbitrary fundamental solution of
$ L$, that is, a solution to the equation   ${
L}^*G(x_0,y_0,x,y)=\delta (x_0,y_0)$. It is well-known that $G$ is
a real-analytic function in $\mathbb{R}^2$ except at the point $P(x_0,
y_0)$. Its continuation to the complex space $\mathbb{C}^2$ has
logarithmic singularities on the complex characteristics passing
through this point, i.e., on $K_P:=\{(x-x_0)^2 + (y-y_0)^2=0\}.$ Our
proof is based on the idea suggested by Garabedian \cite{G} to deform
contour $\gamma$ across the curve $\Gamma$ from the domain $U_1$ to
the domain $U_2$. To be able to realize this deformation, first, we use a special representation for
a fundamental solution, that is,  a sum of two functions, each of those
has a singularity on a single characteristic only. This representation is given by the Theorem \ref{L1} above. Next,
we replace the fundamental solution $G$ with
 a so-called {\it   reflected fundamental solution}. After proving the existence and uniqueness of the reflected fundamental solution, we describe the deformation of $\gamma$ and obtain the desired reflected formula. Finally, we simplify the formula and discuss the cases for which it reduces to the simplest point to point form.

\subsection{The reflected fundamental solution}
\label{sec:3} 

This section is devoted to the construction of  the reflected fundamental solution $\widetilde{G}$, which  plays the key role by enabling us to deform the contour $\gamma$ across the boundary. $\widetilde{G}$ depends on the operator $L$ and the curve $\Gamma$ \footnote{$\widetilde{G}$ depends on the boundary condition as well, but the later is beyond the scope of this paper, see \cite{bs}, \cite{s99} and \cite{bihar} for some relevant results.}. As it will be shown in the next two sections, the reflected fundamental solution determines whether the corresponding reflection formula can be reduced to the point to point form.

Function $\widetilde{G}$ is a solution of the equation ${ L}^* _{\mathbb C}\widetilde{G} =0$
subject to the boundary condition $G=\widetilde{G}$ on
$\Gamma_{\mathbb{C}}$ and has singularities only on the ``reflected'' characteristic
lines $\os{l}_1$ and $\os{l}_2$ (see Fig.~1) intersecting the real space at the reflected point $Q=R(P)$  in the
domain $U_2$ and intersecting $\Gamma_{\mathbb{C}}$ at $K_P\cap
\Gamma_{\mathbb{C}}$.

 We seek the reflected
fundamental solution in the form
\begin{equation}\label{E:1.5}
\widetilde{G} (z_0,\zeta _0,z,\zeta )=
-\frac{1}{4\pi}(\widetilde{G}_1 (z_0,\zeta _0,z,\zeta  )
+\widetilde{G}_2 (z_0,\zeta _0,z,\zeta )),
\end{equation}
where the functions $\widetilde{G}_j$, $j=1,2$ are defined as the
solutions to the following  Cauchy-Goursat problems with prescribed
singularities,
\begin{eqnarray}\label{fundr133}
\begin{cases}
{ L}^* _{\mathbb C}\widetilde{G}_j =0,\qquad j=1,2, \label{ref11} \\
\widetilde{G}_{j_{|_{\Gamma _{\mathbb{C}}}}}= \, G_{j_{|_{\Gamma
_{\mathbb{C}}}}}
, \label{ref12}\\
\widetilde{G}_j\mbox{ has singularities only on the char.
}\;\widetilde{l}_j=:\{ \widetilde{\psi}_j(z,\,\zeta )=0
\},\label{ref13}
\end{cases}
\end{eqnarray}
where $\widetilde{\psi}_1=\widetilde{S}(\zeta )-z_0$ and
$\widetilde{\psi}_2=S(z)-\zeta _0$ are solutions of the
Hamilton-Jacobi equation
\begin{equation}\label{psi}
\pd{\widetilde{\psi} _j}{z}  \cdot \pd{\widetilde{\psi}
_j}{\zeta}=0.
\end{equation}

First, we construct the solutions to the problems (\ref{fundr133}) as some
formal expansions. Then  we justify their convergence.

We seek  this expansion in the form \cite{ludwig},
\begin{equation}\label{Gja}
\widetilde{G}_1=\sum\limits_{k=0}^{\infty}\beta ^1_k(z_0,\zeta
_0,z,\zeta )f_k(\widetilde{\psi}_1),
\end{equation}
\begin{equation}\label{Gjb}
\widetilde{G}_2=\sum\limits_{k=0}^{\infty}\beta ^2_k(z_0,\zeta
_0,z,\zeta )f_k (\widetilde{\psi}_2),
\end{equation}
where $f_k(\xi )$ is defined by \eqref{fk}.

Substituting (\ref{Gja}) and (\ref{Gjb}) into (\ref{ref11}), we
obtain the following recursion for the coefficients $\beta ^1 _k$ and $\beta ^2 _k,$
\begin{equation}\label{transp1}
\begin{aligned}
& \Bigl ( \pd{\beta ^1 _0}{z}-B\, \beta ^1 _0\Bigr
)\widetilde{S}^{\prime}(\zeta )=0,\quad \Bigl (
\pd{}{z}\,\beta ^1 _{k+1}-B\, \beta ^1 _{k+1}\Bigr
)\widetilde{S}^{\prime}(\zeta )=- L^*_{\mathbb
C}\beta ^1 _k, \quad k\ge 0, \\
& \Bigl ( \pd{\beta ^2 _0}{\zeta}-A\, \beta ^2 _0\Bigr
){S}^{\prime}(z )=0,\quad \Bigl ( \pd{}{\zeta}\,\beta ^2
_{k+1}-A\, \beta ^2 _{k+1}\Bigr ){S}^{\prime}(z )=-
L^*_{\mathbb C}\beta ^2 _k, \quad k\ge 0
\end{aligned}
\end{equation}
subject to the following initial conditions
\begin{equation}
\beta ^1 _{k_{|_{\Gamma _{\mathbb{C}}}}}= \,
\alpha^1_{k_{|_{\Gamma _{\mathbb{C}}}}}, \qquad
\beta ^2 _{k_{|_{\Gamma _{\mathbb{C}}}}}= \,
\alpha^2_{k_{|_{\Gamma _{\mathbb{C}}}}},\qquad k=0,1,2,\dots \, .
\end{equation}
Note that both $S^{\prime}(z)$ and $\widetilde{S}^{\prime}(\zeta )$ do
not vanish on $\Gamma _{\mathbb{C}}$, \cite{davis} Chapter 7, p.~42. Indeed, functions 
$S(z)$ and $\widetilde{S}(\zeta )$ are inverse of each other (see (\ref{E:999})), so  $\widetilde{S}(S(z) )=z$.
Differentiating the later equation and taking into account that $S(z)=\zeta$ on $\Gamma _{\mathbb{C}}$,
we obtain $\widetilde{S}^{\prime}(\zeta )\cdot S^{\prime}(z)=1$. In $\mathbb{R}^2$, 
$\widetilde{S}^{\prime}(\zeta )=\overline{S^{\prime}(z)}$, therefore, $|S^{\prime}(z)|=|\widetilde{S}^{\prime}(\zeta )|=1$ on $\Gamma$. Thus, both functions $S^{\prime}(z)$ and $\widetilde{S}^{\prime}(\zeta )$ are nonzero throughout some neighborhood of $\Gamma$ as continuous functions.

Thus, functions
$\beta ^1_k$ and $\beta ^2 _k$ are uniquely determined 
near $\Gamma _{\mathbb{C}}$, specifically
\begin{equation}\label{beta0}
\begin{aligned}
& \beta ^1 _0=\exp \Bigl (\int\limits _{\zeta
_0}^{\zeta}A(\widetilde{S}(\zeta ),\tau )d\tau +\int\limits _{z
_0}^{z}B(t,\zeta  )dt\, +\int\limits _{z _0}^{\widetilde{S}(\zeta
)}[B(t,\zeta _0 )-B(t,\, \zeta )]dt\Bigr ), \\
& \beta ^2 _0=\exp \Bigl (\int\limits _{\zeta
_0}^{\zeta}A(z,\tau )d\tau+ \int\limits _{\zeta
_0}^{S(z)}[A(z_0,\tau )-A(z,\tau )]d\tau  +\int\limits _{z
_0}^{z}B(t,S(z))dt\Bigr ).
\end{aligned}
\end{equation}
Hence, the formal expansions for the functions
$\widetilde{G}_1,\,\widetilde{G}_2$ satisfying conditions
(\ref{ref11})  are constructed.

\begin{lemma}\label{lemma}  The series \eqref{Gja} and \eqref{Gjb} converge near $\Gamma
_{\mathbb{C}}$.
\end{lemma}
{\bf Proof of the Lemma \ref{lemma}}

Let us prove the convergence of the series \eqref{Gjb} by
considering an auxiliary family of problems depending on parameter
$\xi$:
\begin{equation}\label{auxip}
\begin{cases}
& { L}^* _{\mathbb C}V_{\xi}(z_0,\zeta _0,z,\zeta ,\xi ) =0,  \\
& {V_{\xi}}(z_0,\zeta _0,z,S(z) ,\xi ) = \Phi (z_0,\zeta _0,z,S(z),\xi )
, \\
& {V_{\xi}}(z_0,\zeta _0,\widetilde{S}(\zeta _0 -\xi),\zeta ,\xi )=0
.
\end{cases}
\end{equation}
Here $\Phi$ is a given analytic function, that has Taylor expansion
$$
\Phi (z_0,\zeta _0,z,S(z),\xi ) =\sum\limits _{k=0}^{\infty}
\alpha ^2 _k (z_0,\zeta _0, z, S(z) ) \frac{(S(z)-\zeta _0+\xi
)^{k+1}}{(k+1)!},
$$
where coefficients $\alpha _k ^2(z_0,\zeta _0, z, \zeta )$ are the same as in \eqref{anal} \cite{fundsol}.

Taylor expansion of the solution to the problem \eqref{auxip}
(if it exists) has the form
\begin{equation}\label{auxis}
V_{\xi}(z_0,\zeta _0,z,\zeta ,\xi ) =\sum\limits _{k=0}^{\infty}
\beta ^2 _k (z_0,\zeta _0, z, \zeta ) \frac{(S(z)-\zeta _0+\xi
)^{k+1}}{(k+1)!}
\end{equation}
where the coefficients $\beta ^2 _k$ are the same as  the
coefficients in  series \eqref{Gjb}. Convergence of the later,
therefore, followed from convergence \eqref{auxis}. To show
existence and uniqueness of the solution to the problem \eqref{auxip} in
the class of analytic functions we use the 
substitution
\begin{equation}\label{density}
V_{\xi}(z_0,\zeta _0,z,\zeta ,\xi ) =\int\limits _{S(z)
}^{\zeta}d\tau\int\limits _{\widetilde{S}(\zeta _0 -\xi) }^{z }\mu
(z_0,\zeta _0,t,\tau ,\xi )dt+\Phi (z_0,\zeta _0,z,S(z),\xi )
\end{equation}
with unknown density $\mu$, which reduces the problem \eqref{auxip}
to the
Volterra integral equation
\begin{equation}\label{inteq1}
\begin{aligned}
& \mu (z_0,\zeta _0,z,\zeta ,\xi )+A(z,\zeta
)S^{\prime}(z)\int\limits _{\widetilde{S}(\zeta _0-\xi )
}^{z}\mu (z_0 ,\zeta _0, t,S(z) ,\xi )dt\\
& -A(z,\zeta )\int\limits _{S(z) }^{\zeta }\mu (z_0 ,\zeta _0,
z,\tau ,\xi )d\tau -B(z,\zeta )\int\limits _{\widetilde {S}(\zeta
_0 -\xi ) }^{z }\mu (z_0 ,\zeta _0, t,\zeta ,\xi )dt \\
& -F(z,\zeta )\int\limits _{S(z) }^{\zeta }d\tau \int\limits
_{\widetilde {S}(\zeta _0 -\xi ) }^{z }\mu (z_0 ,\zeta _0, t,\tau
,\xi )d\tau=\Psi (z_0,\zeta _0,z,\zeta ,\xi ),
\end{aligned}
\end{equation}
where
\begin{eqnarray}
F(z,\zeta )=\pd{}{z}A(z,\zeta)+\pd{}{\zeta}B(z,\zeta )-C(z,\zeta
),\\
\Psi (z_0,\zeta _0,z,\zeta ,\xi )=F(z,\zeta )\Phi (z_0,\zeta
_0,z,\xi )+A(z,\zeta )\pd{}{z}\Phi (z_0,\zeta _0,z,\xi ).
\end{eqnarray}
The existence and uniqueness of the analytic solution of  equation \eqref{inteq1} can be  proven by iteration technique
described in \cite{vekua}, Chapter 1, p.~11. Thus, there exists unique solution of
\eqref{auxip}, which has unique Taylor expansion with respect to
variable $\xi$ at the point $\xi =-(S(z)-\zeta _0)$, this expansion coincide with the
expansion \eqref{auxis}. Thus, series \eqref{auxis} converges in the
neighborhood of $\Gamma$, and so does \eqref{Gjb}.

Analogously, considering the following auxiliary problem depending on
parameter $\eta$:
\begin{equation}\label{auxip2}
\begin{cases}
& { L}^* _{\mathbb C}V_{\eta}(z_0,\zeta _0,z,\zeta ,\eta ) =0,  \\
& {V_{\eta}}(z_0,\zeta _0,\widetilde{S}(\zeta),\zeta,\eta)=
\sum\limits _{k=0}^{\infty}\alpha _k^1(z_0,\zeta
_0,\widetilde{S}(\zeta),\zeta)
 \frac{(\widetilde{S}(\zeta ) -z _0 +\eta  )^{k+1}}{(k+1)!} , \\
& {V_{\eta}}(z_0,\zeta _0,z,S(z _0-\eta) ) ,\eta ) =0,
\end{cases}
\end{equation}
whose solution has  Taylor expansion $V_{\eta}=\sum
_{k=0}^{\infty}\beta ^1 _k \frac{(\widetilde{S}(\zeta ) -z _0
+\eta )^{k+1}}{(k+1)!}$, one can show  convergence of \eqref{Gja}.
 That finishes the proof.


\subsection{The reflected fundamental solution as a multiple-valued function}
\label{sec:33} 

As it was conjectured in \cite{alice}: ``Perhaps Looking-glass milk isn't good to drink''.  In this section we show that the reflected fundamental solution (the looking-glass fundamental solution), except for some special cases, does not inherit all of the properties of  a ``true'' fundamental solution, in particular, the representation (\ref{r1}) with the Riemann function as a factor of the logarithm
does not hold. Moreover, as we are about to show,   the factors of the logarithms in $\widetilde{G}_1$ and $\widetilde{G}_2$ (see (\ref{Gja}) and (\ref{Gjb})) are not the same, which makes the reflected fundamental solution a multiple-valued function even in $\mathbb{R}^2$. 
The later  explains (see Section \ref{sec:34}) why point to point reflection almost always fails.

Indeed, consider a point moving along a continuous curve $\gamma$ surrounding either the branch line $z=z_0$ of $G_1$ or  the branch line $\zeta =\zeta _0$ of $G_2$ (\ref{G}). 
As the point makes a complete cycle around the line, it passes to the next sheet of  Riemann surface,  while  going around the  cyclic path surrounding both characteristics at once, it remains on the same sheet of the Riemann surface.
For the reflected fundamental solution $\widetilde{G}$, the point passes to the next sheet even if the curve $\gamma$ lays in $\mathbb{R}^2$ and surrounds both intersecting branch lines,  $\widetilde{S}(\zeta )=z_0$ and
$S(z)=\zeta _0$.

To show this, let us compute the increment of the function $\widetilde{G}$ when a curve $\gamma \subset\mathbb{R}^2$ is a circle of a small radius $\rho$ centered at the point $Q \bigl ( \widetilde{S}(\bar{z}_0), S(z_0) \bigr )$:
\begin{equation}\label{increment}
2\pi i (-\frac{1}{4\pi})\sum\limits _{j=0}^{\infty}\Bigl ( \beta _j^1 \frac{(\widetilde{S}(\zeta )-z_0)^j}{j!}-
\beta _j^2 \frac{({S}(z )-\zeta _0)^j}{j!} \Bigr ) .
\end{equation}
Taking into account that $\zeta = \bar{z}$ in $\mathbb{R}^2$,  let us set $z=\widetilde{S}(\bar{z}_0)+\rho e^{i\phi}$ and 
$\zeta={S}({z}_0)+\rho e^{-i\phi}$, and expand the Shwarz function and its inverse into Taylor series at the point $Q$: $S(z)=\bar{z}_0+C_1\rho e^{i\phi} + o(\rho )$, $\widetilde{S}(\bar{z})={z}_0+\bar{C}_1\rho e^{-i\phi} + o(\rho )$.

Without loss of generality assume that the coefficients $A$, $B$ and $C$ in (\ref{eqc}) are constants (otherwise we should use their Taylor expansions in this analyses), then
$\beta _0 ^1=\beta _0 ^2$ (see (\ref{beta0})), and 
\begin{eqnarray}
&\beta ^1 _1 =(AB-C)e^{(A(\zeta -\zeta _0)+B(z-z_0))}\bigr ( (z-\widetilde{S}(\zeta ))/\widetilde{S}^{\prime}(\zeta )+\zeta -\zeta _0\bigl ), \\
&\beta ^1 _2 =(AB-C)e^{(A(\zeta -\zeta _0)+B(z-z_0))}\bigr ( (\zeta-\widetilde{S}(z ))/{S}^{\prime}(z )+z -z _0\bigl ).
\end{eqnarray}
Thus, the increment (\ref{increment}) becomes
\begin{equation}\label{increment1}
\begin{aligned}
\frac{\rho}{2i}(AB-C)\Bigl ( [C_1\widetilde{S}(\bar{z}_0)-C_1z_0 +S(z_0)-\bar{z}_0  ]\bar{C}_1 e^{-i\phi}\\
- [\bar{C}_1{S}({z}_0)-\bar{C}_1\bar{z}_0 +\widetilde{S}(\bar{z}_0)-{z}_0]{C}_1 e^{i\phi}\Bigr )+ o(\rho ).
\end{aligned}
\end{equation}
Formula (\ref{increment1}) shows that the increment can be equal zero only in two cases: either when (i)  $AB-C=0$ or (ii)
expressions in the brackets equal zero. The later happens if boundary $\Gamma$ is a segment of a straight line, while (i), for example, holds if operator $L$ is the Laplacian. 

Having the detailed description of the reflected fundamental solution we are ready to derive the reflection formula by explaining how the contour $\gamma$ in (\ref{E:1.20}) can be deformed from one side of the reflecting surface $\Gamma _{\mathbb{C}}$ to the other.

\subsection{Deformation of the contour}\label{sec:34}

Formula (\ref{E:1.20}) involves integration over a contour $\gamma\subset U_1$ surrounding both  characteristics on which functions $G_1$ and $G_2$ have singularities (lines $l_1$ and $l_2$ Fig.~1).
\begin{figure}[h]
\setlength{\unitlength}{1cm}
\begin{minipage}[t]{6.0cm}
\scalebox{.5}[.5] {\includegraphics{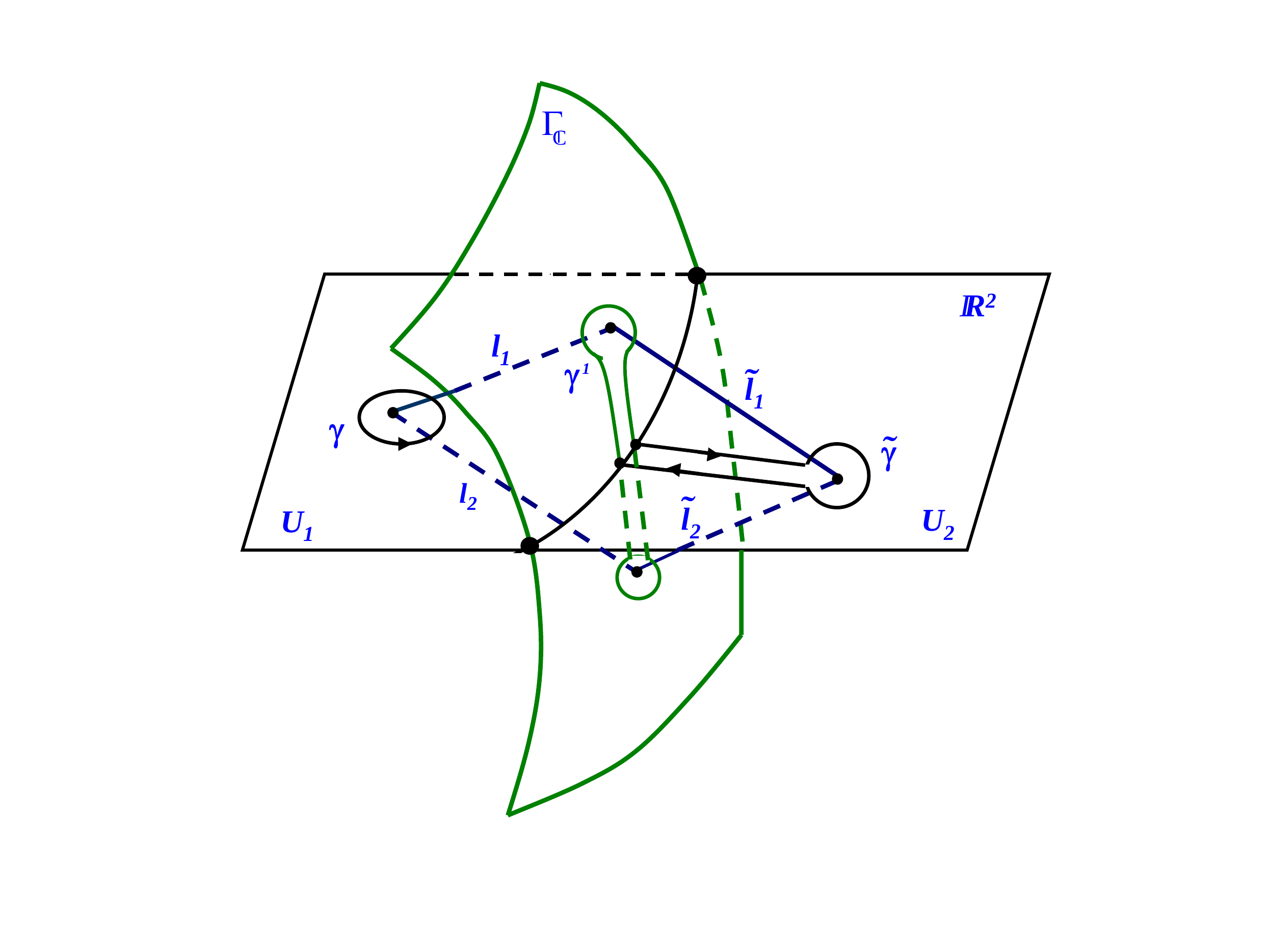}} \caption{Contour
deformation}\label{fig:contourdeformation}
\end{minipage}
\end{figure}
To express value $u(P)$ in terms of values of $u(x,y)$ in $U_2$, that is,
to construct a reflection formula,  
it is sufficient to deform the contour $\gamma$ from the domain $U_1$ to the domain
$U_2$. Note that since the integrand in (\ref{E:1.20}) is a
closed form, $d\omega =0$, the value of the integral will not change while we are
deforming the contour $\gamma$ homotopically. 

First, the contour is  deformed to the
complexified curve $\Gamma_{\mathbb{C}}$. 
Taking into account that the characteristics of $G$ passing through the point $P$ intersect $\Gamma _{\mathbb{C}}$ at two different points 
in ${\mathbb{C}}^2$,
assume that the point $P$ lies so close to the curve $\Gamma$ that
there exists a connected, univalently projected onto a  plane,
domain $\Omega\subset \Gamma _{\mathbb{C}}$ that 
contains both points of intersections
\cite{sss}. 

We start the deformation with stretching the contour $\gamma$ (see (\ref{E:1.20})) in the real plane until its  small arc  reaches the curve $\Gamma$
(it becomes a mirror image of $\os{\gamma}$ in Fig. 1).
Then we substitute  a sum of $G_1$ and $G_2$ for $G$ in (\ref{E:1.20}) and split the integral:
\begin{equation}\label{E:1.20split}
    u(P)
=\int\limits_{\gamma} \omega [ u,G_1  ]+\int\limits_{\gamma} \omega [ u,G_2  ].
\end{equation}
Note that contour $\gamma$ is not closed on  Riemann surfaces of each $G_1$ and $G_2$ (see Section \ref{sec:33}).  As a point of disconnection (one of two endpoints) let us choose a point $K\in \gamma\cap \Gamma$. Then in the first integral in (\ref{E:1.20split}) we ``lift'' the contour $\gamma$ to  $\gamma ^{\prime}$ (solid line above the plane  in Fig. 1) such that we do not move some points of $\gamma\cap \Gamma$ in the neighborhood of the point $K$.  Then we do  the symmetric (with respect to plane $\mathbb{R}^2$) deformation in the second integral in (\ref{E:1.20split}).

Taking into account that $u_{\mid _{ \Gamma _{ \mathbb{C} } } }=0$, differential form $\omega$
   \eqref{E:1.20b} on $\Gamma_{\mathbb{C}}$ becomes
\begin{equation}\label{E:1.20bg}
   \omega ^{\prime} [  u,G_j ]=  G_j\pd{u}{y}\, dx
      - G_j\pd{u}{x}\, dy, \quad j=1,2.
\end{equation}


Now we can replace $G_j$ with  $\widetilde{G}_j$ (see formula  (\ref{E:1.5})).
Indeed,  according to (\ref{fundr133})
\begin{equation}
\int\limits_{\gamma ^{\prime}}\omega ^{\prime} [  u,G_1 ] =\int\limits_{\gamma ^{\prime}}\omega ^{\prime} [  u,\widetilde{G}_j ].
\end{equation} 
In order to deform contour $\gamma ^{\prime}$ from $\Gamma
_{\mathbb{C}}$ to the domain $U_2$, it is necessary to apply the
``mirror'' deformation procedure. Note that during this deformation the point $K$ is fixed and contour surrounds one of the ``reflected'' characteristic
lines $\os{l}_1$ or $\os{l}_2$  (see Fig.~1) intersecting the real space at the reflected point $Q=R(P)$  in the
domain $U_2$ and intersecting $\Gamma_{\mathbb{C}}$ at $K_P\cap
\Gamma_{\mathbb{C}}$.



Finally, we have
\begin{equation}\label{E:1.20splitr}
    u(P)
=\int\limits_{\os{\gamma}} \omega [ u,\widetilde{G}_1  ]+\int\limits_{\os{\gamma}} \omega [ u,\widetilde{G}_2  ].
\end{equation}
This formula can be rewritten as a single integral
\begin{equation}\label{E:1.20splitr2}
    u(P)
=\int\limits_{\os{\gamma}} \omega [ u,\widetilde{G}  ],
\end{equation}
but as it was discussed in Section \ref{sec:33} the contour $\os{\gamma}$, generally, is not closed on Riemann surface of 
$\widetilde{G}$, so in most of the cases we do not expect to be able to move the point $K$ (see Fig. 2) from the curve $\Gamma$.  Formula (\ref{E:1.20splitr2}) is a version of a desired reflection formula. In the next section we simplify it
and show that it holds in the large.


\subsection{The reflection formula in the large}


Formula (\ref{E:1.20splitr2}) in $(z,\zeta )$ variables has the form
\begin{equation}\label{E:1.20d}
    u(P)
=\int\limits_{\os{\gamma}} \os{\omega}[u,\os{G}], 
\end{equation}
where
\begin{equation}\label{E:1.20dd}
  \os{\omega}[u,\os{G}]   
=i \Bigl ( \bigl \{ u\pd{\os{G}}{\zeta}-\os{G}\pd{u}{\zeta}-2Au\os{G}
\bigr \} d\zeta
      - \bigl \{u\pd{\os{G}}{z}-\os{G}\pd{u}{z}-2Bu\os{G} \bigr \} dz
\Bigr ) .
\end{equation}
Here  $\os{G}$ is the reflected fundamental solution and
contour $\os{\gamma}\subset U_2$ surrounds the point $Q$ (see Fig.~2).
\begin{figure}[h]
\setlength{\unitlength}{1cm}
\begin{minipage}[t]{6.0cm}
\scalebox{.4}[.4] {\includegraphics{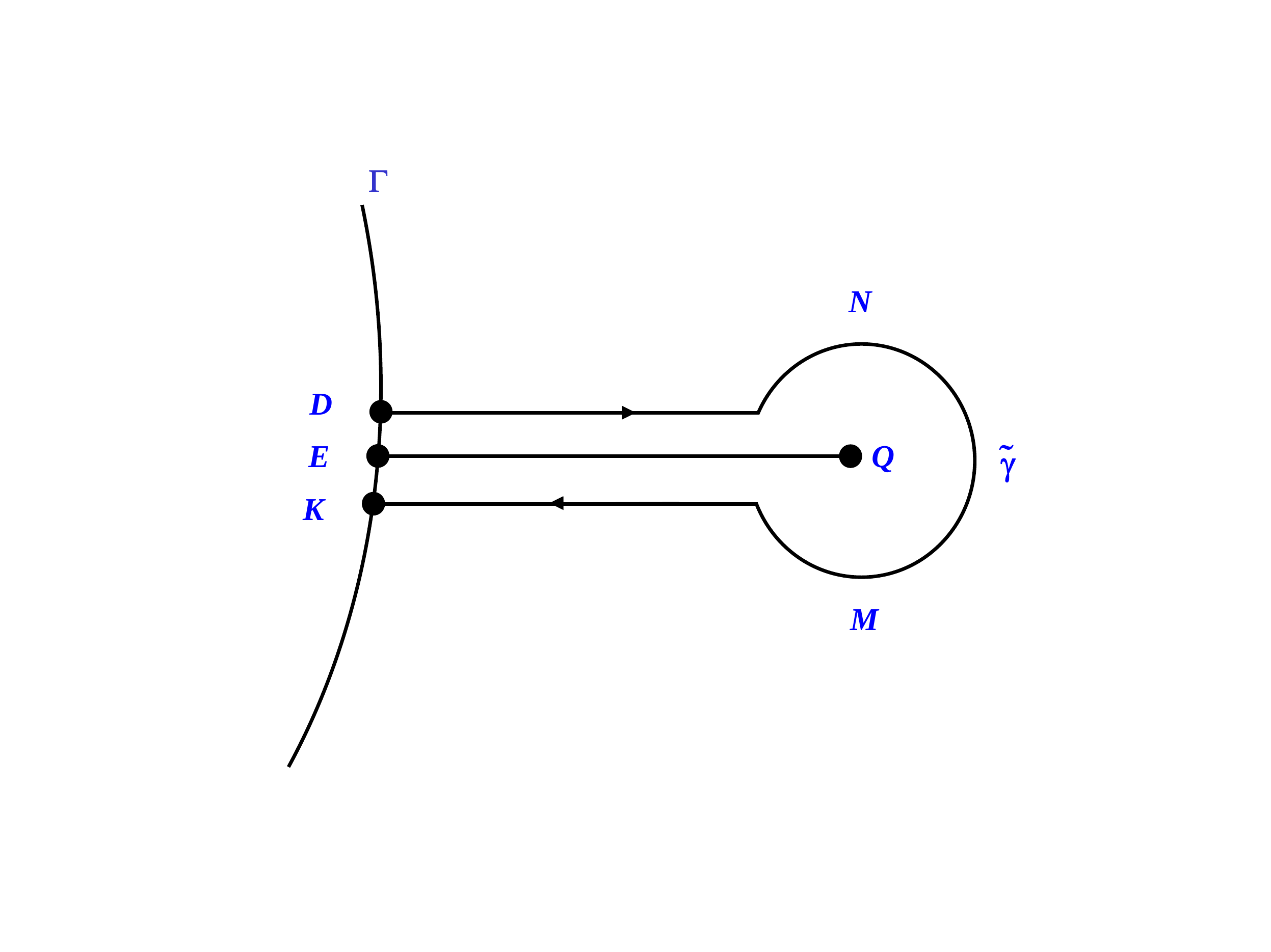}} \caption{Contour
transformation}\label{fig:contourtransformation}
\end{minipage}
\end{figure}
Recall that $\os{G}$ is a sum of two series (with certain radii of convergence). 
Now we are going to show that the formula holds in the large. 

Let us rewrite functions $\os{G}_l$ in the form:
\begin{equation}\label{E:1.21d}
\os{G}_l=V_l\ln \os{\psi}_l+\os{V}_l, \qquad l=1,2,
\end{equation}
where
\begin{equation}\label{E:1.22d}
V_l=\sum\limits _{j=0}^{\infty}\beta ^l _j\frac{(\os{\psi}_l)^j}{j!}, \qquad
\os{V}_l=\sum\limits _{j=0}^{\infty}\beta ^l _j\frac{(\os{\psi}_l)^j}{j!}C_j,
\end{equation}
and
\begin{equation}\label{psisi}
\widetilde{\psi}_1=\widetilde{S}(\zeta )-z_0, \qquad
\widetilde{\psi}_2=S(z)-\zeta _0.
\end{equation}
Substituting  (\ref{E:1.21d}) and  (\ref{E:1.22d}) into (\ref{E:1.20d}) and letting the radius of the arc NM to zero (see Fig.~2)  result in vanishing integrals of the terms involving products of function $\os{V}_l$ and derivatives of function $u$  as integrals of  holomorphic functions over a closed contour. Combining terms in (\ref{E:1.20d})  involving derivatives of  logarithms and separating them from the terms  involving logarithmic functions yields
\begin{equation}\label{QI}
u(P)={\mathbb{Q}}+{\mathbb{I}},
\end{equation}
where ${\mathbb{Q}}$ and ${\mathbb{I}}$ in characteristic variables have the form
\begin{equation}\label{E:1.25d}
    {\mathbb{Q}}
=-\frac{i}{4\pi}\sum\limits _l\int\limits_{\os{\gamma}} \Bigl ( 
uV_l\pd{}{\zeta}\bigl( \ln {\os{\psi}_l}\bigr )d\zeta -
  uV_l\pd{}{z}\bigl ( \ln {\os{\psi}_l} \bigr )dz 
 \Bigr ),
\end{equation}

\begin{equation}\label{E:1.26d}
    {\mathbb{I}}
=-\frac{1}{4\pi}\int\limits_{\os{\gamma}}  \os{\omega}[u,\os{G}_1]\ln {\os{\psi}_1}\,-\frac{1}{4\pi}\int\limits_{\os{\gamma}}  \os{\omega}[u,\os{G}_2]\ln {\os{\psi}_2}\,.
\end{equation}
Substituting series \eqref{E:1.22d} for $V_l$ into \eqref{E:1.25d} and computing the residues at the point $Q$ where the integrand has the simple pole, we have
\begin{equation}\label{E:1.27d}
\begin{aligned}
    &{\mathbb{Q}}
=-\frac{1}{2}u(Q)\Bigl ( \exp \bigl ( \int\limits_{z_0}^{\os{S}(\zeta _0)} 
B(t,S(z_0))dt + \int\limits_{\zeta _0}^{{S}(z _0)} 
A(z_0,\tau )d\tau \bigr )\\
& +\exp \bigl ( \int\limits_{\zeta _0}^{{S}(z _0)} 
A({\os{S}(\zeta _0)},\tau)d\tau + \int\limits_{z _0}^{\os{S}(\zeta _0)} 
B(t,\zeta _0 )dt \bigr )   
 \Bigr ),
\end{aligned}
\end{equation}
which holds in the large.

Using properties of the logarithmic function and replacing the contour $\os{\gamma}$ with a segment EQ, second integral can be rewritten as 
\begin{equation}\label{E:1.28d}
    {\mathbb{I}}
=2\pi i (-\frac{1}{4\pi})  \Bigl (  
\int\limits_{E}^Q   \os{\omega}[u,V_1]-\int\limits_{E}^Q   \os{\omega}[u,V_2]\Bigr ) 
=\frac{1}{2i}\int\limits_{E}^Q   \os{\omega}[u,V],
\end{equation}
where $V=V_1-V_2$. Note, that  the logarithms in (\ref{E:1.26d}) have complex conjugated arguments
(\ref{psisi}) in $\mathbb{R}^2$, however they  cancel each other only if the factors $V_1$ and $V_2$ are equal, which generally is not the case.

Even though the later formula involves series $V_1$ and $V_2$, it is also holds in the large, since these expansions can be interpreted as solutions of the following Cauchy problems
\begin{eqnarray}\label{glob}
\begin{cases}
{ L}^* _{\mathbb C}{V}_j =0,\qquad j=1,2,  \\
V_{j_{|_{\Gamma _{\mathbb{C}}}}}= \, {\mathfrak R}_{{|_{\Gamma
_{\mathbb{C}}}}}
, \\
V_j=\exp \Bigl \{\int\limits _{\zeta _0} ^{
\zeta} A(\theta _j,\tau )d\tau +
\int\limits _{z _0} ^{z} B(t,\eta _j )dt\Bigr \}\mbox{  on the characteristic
}\;\widetilde{l}_j=:\{ \widetilde{\psi}_j(z,\,\zeta )=0
\},
\end{cases}
\end{eqnarray}
where $\theta _1 ={\os S}(\zeta)$, $\theta _2 =z$, $\eta _1=\zeta$ and $\eta _2 =S(z)$. 
Problem \eqref{glob} by a substitution with unknown density $\mu $, for example, for $j=2$ 
\begin{equation}\label{density2}
\begin{aligned}
& V_{2}(z_0,\zeta _0,z,\zeta ) =\int\limits _{S(z)
}^{\zeta}d\tau\int\limits _{\widetilde{S}(\zeta _0 ) }^{z }\mu
(z_0,\zeta _0,t,\tau )dt \\
& +{\mathfrak R}(z_0,\zeta _0,z,S(z))\, 
e^{ \int\limits _{\zeta _0}^{\zeta } A(z,\tau )\, d\tau \, -  \int\limits _{\zeta _0}^{{S}(z)} A(z,\tau)\, d\tau  },
\end{aligned}
\end{equation}
can be reduced to the Volterra integral equation,  
whose solution as a function of four complex variables exists and unique in some cylindrical domain near $\Gamma$ (see \cite{vekua} Chapter 1, p.~11).
Thus, the solutions of \eqref{glob} exist in $\mathbb{C}^4$ as  multiple-valued analytic functions, whose singularities coincide with those of $S(z)$ and $\os{S}(\zeta )$.

Combining \eqref{QI}, \eqref{E:1.27d} and \eqref{E:1.28d}
we arrive at the formula \eqref{formula}, which  proves the theorem.

\section{Conclusions and remarks}\label{sec5}
\subsection{Equations with constant coefficients}

We have obtained a reflection formula for elliptic equations with analytic coefficients subject to homogeneous Dirichlet conditions on a real analytic curve. This is a point to compact set reflection, which in some cases can be essentially simplified.

Consider the case when the coefficients $a$, $b$ and $c$ in the equation \eqref{intro1} are constants,
\begin{equation}\label{concl1}
\Delta _{x,y}u+a\pd{u}{x}+b\pd{u}{y}+cu=0,
\end{equation}
 and, therefore, $A$, $B$ and $C$ are constants as well. In this case solutions $\alpha _k^j$ to the problems \eqref{transp} -- \eqref{initcond} can be written explicitly, and the Riemann function \eqref{anal}  has the form  
\begin{equation}\label{analConst}
{\mathfrak R}(z_0 ,\zeta _0 , z, \zeta )=\sum\limits _{k=0}^{\infty}
 \frac{\bigl ( (z-z_0) (\zeta-\zeta _0)(AB-C)\bigr )^k}{(k!)^2}\exp{(A(\zeta-\zeta _0)+B(z-z_0))}.
\end{equation}

Our main conclusion confirms the fact that the point to point reflection is quite rare.
\begin{theorem}\label{negative}
For non-trivial solutions of elliptic equation \eqref{concl1}  with constant coefficients vanishing on a real-analytic curve $\Gamma$, there is no point to point reflection unless one of the following conditions hold:

(i) $\Gamma$ is a line,

(ii) $a^2+b^2-4c=0$.
\end{theorem}
{\bf Proof: } The proof immediately follows from the fact that the integral term ${\mathbb{I}}\ne 0$ in  \eqref{E:1.28d}. Indeed, formula (\ref{increment1}) imply that $V\ne 0$. Thus, for  ${\mathbb{I}}$ to be zero, function $u$ and its first derivative must vanish on a path joining the curve $\Gamma$ with the reflected point, which contradicts the assumption that $u$ is not equal zero identically.

\begin{theorem}\label{linia}
Let   $\Gamma := \{ \alpha x+\beta y + \delta=0 \}$ be a line. Then for any solution of the equation
$\Delta u +au_x+bu_y+cu=0$ with constant coefficients vanishing on $\Gamma$ the following point to point reflection formula holds in $\mathbb{R}^2$:
\begin{equation}\label{linlin}
u(P)=-\exp {\Bigl (-\frac{(\alpha x_0+\beta y_0 + \delta)(a\alpha +b\beta )}{\alpha ^2+\beta ^2}\Bigr )}u(Q).
\end{equation}
\end{theorem}

{\bf Proof: } Under the assumptions of the theorem, the Shwarz function is $S(z)=mz+q$, 
where 
$$
m=\frac{\beta ^2 -\alpha ^2+i\, 2\alpha\beta}{\alpha ^2 +\beta ^2}, \qquad
q=\frac{ -2\alpha \delta +i\, 2\beta\delta}{\alpha ^2 +\beta ^2}.
$$
Functions $V_1$ and $V_2$ are equal (see \eqref{glob}),
\begin{equation}\nonumber
V_1=V_2 =\sum\limits _{k=0}^{\infty}
 \frac{\bigl ( (mz+q-\zeta _0)({\bar m}\zeta +{\bar q} -z_0) (AB-C)\bigr )^k}{(k!)^2}\, e^{(A(\zeta-\zeta _0)+B(z-z_0))},
\end{equation}
and
therefore, $V=V_1-V_2=0$, and the integral ${\mathbb{I}}=0$ (see \eqref{E:1.28d}).
 Formula \eqref{E:1.27d} can be simplified, and
${\mathbb{Q}}=-u(Q)\, e^{ A(mz_0+q-\zeta _0)+B({\bar m}\zeta +{\bar q} -z_0) }$. The later in variables $(x,y)$ 
gives \eqref{linlin}.

\begin{corollary}
Let $\Gamma$ be a  line with equation $y=0$, then for any solution of  \eqref{concl1}  vanishing on $\Gamma$ the following reflection formula holds
\begin{equation}
u(x_0,y_0)= - e^{-by_0}u(x_0,-y_0).
\end{equation}
\end{corollary}
\begin{corollary}
Let $\Gamma$ be a  line with equation $x=0$, then for any solution of  \eqref{concl1}  vanishing on $\Gamma$ reflection formula has the form
\begin{equation}
u(x_0,y_0)= - e^{-ax_0}u(-x_0,y_0).
\end{equation}
\end{corollary}
\begin{corollary}
If $a=b=0$ formula \eqref{linlin}  recovers known point to point reflection for solutions of the Helmholtz equation vanishing on a  line
 $$u(P)=-u(Q).$$
\end{corollary}
\begin{remark} 
Note, that
in the case of the Helmholtz equation, $a=b=0$ and $c=\lambda ^2$, when $\Gamma$ is a real-analytic curve, formula  \eqref{E:1.27d} can be simplified, and
${\mathbb{Q}}=-u(Q)$, but ${\mathbb{I}}\ne 0$ in  \eqref{E:1.28d} unless $\Gamma$ is a  line \cite{laguna} Chapter 9, p.59;
\cite{khavinsons}; \cite{sss}.
\end{remark}

\begin{theorem}
Let   $\Gamma $ be a real-analytic curve. Then for any solution of the equation
$\Delta u +au_x+bu_y+{(a^2+b^2)}/4\, u=0$ vanishing on $\Gamma$ the following point to point reflection formula holds in $\mathbb{R}^2$:
\begin{equation}\label{ptpnew}
u(P)=-e^{A(S(z_0)-\zeta _0)+B({\os S}(\zeta _0)-z_0)}u(Q).
\end{equation}
\end{theorem}
{\bf Proof: } In characteristic variables condition $c={(a^2+b^2)}/4$ is equivalent to $AB-C=0$. Then
 the Riemann function \eqref{analConst} has the simplest form
\begin{equation}
{\mathfrak R}(z_0 ,\zeta _0 , z, \zeta )=e^{A(\zeta-\zeta _0)+B(z-z_0)},
\end{equation}
 and $V_1=V_2={\mathfrak R}$ for any analytic curve $\Gamma$. Thus, the reflection formula has the point to point form \eqref{ptpnew}. 
\begin{remark}
Equation $\Delta u +au_x+bu_y+{(a^2+b^2)}/4\, u=0$ can be transformed into the Laplace equation using the substitution $u(x,y)=v(x,y)e^{-(ax+by)/2}$, where $v$ is a harmonic function, and, therefore, $v$ enjoys the celebrated Schwarz symmetry principle \eqref{E:0.1}. 
\end{remark}
\begin{example}
Formula \eqref{ptpnew} for the unit circle centered at the origin can be rewritten in $(x,y)$ variables as follows
\begin{equation}\label{ptpcircle}
u(x_0,y_0)=-\exp{\Bigl ( \frac{2(ax_0+by_0)(1-x_0^2-y_0^2)}{x_0^2+y_0^2}\Bigr )}u(\frac{x_0}{x_0^2+y_0^2},\frac{y_0}{x_0^2+y_0^2}).
\end{equation}
\end{example}

\subsection{A final remark}

Thus, for elliptic equations of the second order with real-analytic coefficients in $\mathbb{R}^2$, there is no point to point reflection with respect to a real-analytic curve  $\Gamma$ unless $\Gamma$ is a line or the following constrain $a^2+b^2-4c=0$ for the coefficients of the equation holds. 

As it  follows from \cite{lopez} for elliptic equations in 
$\mathbb{R}^2$, there is no point to finite set reflection as well.

Point to compact set reflection is always possible. This set is a curve having  one of its endpoints on a reflecting curve. The other endpoint is located at the reflected point itself.

\textbf{Acknowledgments.} Research of the author was supported in part by OU
Research Challenge Program, award \# RC-09043. The author is especially grateful to the anonymous referee, whose comments  have improved the paper.

\providecommand{\bysame}{\leavevmode\hbox to3em{\hrulefill}\thinspace}


\begin{thebibliography}{1}
\bibitem{as}
D.~Aberra, T.~Savina, \emph{The Schwarz reflection principle for
polyharmonic functions in $\mathbb{R}^2$}, Complex Var. Theory
Appl., \textbf{41} (2000), no 1, 27-44.
\bibitem{bs}
B.P.~Belinskiy and T.V.~Savina, \emph{The Schwarz reflection principle for harmonic functions in $\mathbb{R}^2$ subject to the Robin condition}, J. Math. Anal.  Appl., \textbf{348} (2008), 685-691.
\bibitem{bramble}
J.~Bramble,  \emph{Continuation of biharmonic functions across circular
arcs}, J. Math. Mech., \textbf{7} (1958), N 6, 905--924.
\bibitem{alice}
L.~Carroll, \emph{Through the Looking-glass}, in: Alice in wonderland,
Wordsworth Edition,  1995.
\bibitem{colton}
D.~Colton and R.P.~Gilbert,  \emph{Singularities of solutions to
elliptic partial differential equations with analytic coefficients},
Q. J. Math.  \textbf{19} 1 (1968), 391-396.
\bibitem{davis}
Ph.~Davis, \emph{The Schwarz function and its applications},
Carus Mathematical Monographs, MAA, 1979.
\bibitem{duffin}
R.J.~Duffin, \emph{Continuation of biharmonic functions by reflection},
Duke Math. J., \textbf{22} (1955), N 2, 313--324.
\bibitem{ek}
P.~Ebenfelt and D.~Khavinson,  \emph{On point to point reflection of
harmonic functions across real analytic hypersurfaces in $\mathbb{R}^n$},
J. d\'{}Analyse Math\'{e}matique, \textbf{68} (1996), 145--182.
\bibitem{ebenfelt}
P.~Ebenfelt, \emph{Holomorphic extension of solutions of elliptic partial
differential equations and a complex Huygens principle}, J. London Math.
Soc., \textbf{55} (1997), 87-104.
\bibitem{farwig}
R.~Farwig, \emph{A note on a reflection principle for the biharmonic equation and the Stokes system},
Acta Appl. Math., \textbf{37} (1994), 41--51.
\bibitem{G}
P.R.~Garabedian,  \emph{Partial differential equations with more than two
independent variables in the complex domain}, J. Math. Mech., \textbf{9}
(1960), 241--271.
\bibitem{Gar}
P.R.~Garabedian,  \emph{Partial differential equations}, John Wiley and Sons, Inc., 1964.
\bibitem{hadamard}
J.~Hadamard,  \emph{Lectures on Cauchy's problem in linear partial differential equations},
 Yale University Press, New Haven, 1923.
\bibitem{john}
F.~John,   \emph{Continuation and reflection of solutions of partial
differential equations}, Bull. Amer. Math. Soc., \textbf{63} (1957),
327--344.
\bibitem{john1}
F.~John,   \emph{Plane waves and spherical means applied to  partial
differential equations}, Springer-Verlag, New York-Berlin, 1981.
\bibitem{john2}
F.~John,   \emph{The fundamental solution of linear elliptic
differential equations with analytic coefficients}, Comm. Pure and
Appl. Math., \textbf{2} (1950), 213--304.
\bibitem{laguna}
D.~Khavinson, \emph{Holomorphic partial differential equations and classical potential theory},
Universidad de La Laguna, 1996.
\bibitem{khavinsons}
D.~Khavinson and H.S.~Shapiro,  \emph{Remarks on the reflection
principles for harmonic functions}, J. d\'{}Analyse
Math\'{e}matique,
\textbf{54} (1991),  60--76.
\bibitem{lewi}
H.~Lewi,  \emph{On the reflection laws of second order differential
equations in two independent variables}, Bull. Amer. Math. Soc.,
 \textbf{65} (1959), 37--58.
\bibitem{lopez}
R.R.~L\'opez, \emph{On reflection principles supported on a final set}, J. Math. Anal. Appl., \textbf{351} (2009), 556-566.
\bibitem{ludwig}
D.~Ludwig, \emph{Exact and Asymptotic solutions of the Cauchy
problem}, Comm. Pure Appl. Math., \textbf{13} 3, (1960), 473--508.
\bibitem{poritsky}
H.~Poritsky, \emph{Application of analytic functions to two-dimensional
biharmonic analysis}, Trans. Amer. Math. Soc.,
\textbf{59} (1946), N 2, 248--279.







\bibitem{sss}
T.V.~Savina, B.Yu.~Sternin  and V.E.~Shatalov,  \emph{On a reflection
formula for the Helmholtz equation}, J.  Comm. Techn.  Electronics, \textbf{38}
(1993), no. 7, 132--143.
\bibitem{sss1}
T.V.~Savina, B.Yu.~Sternin  and V.E.~Shatalov,  \emph{On the reflection
law for the Helmholtz equation}, Dokl. Math., \textbf{45}
(1992), no. 1, 42--45.
\bibitem{s99}
T.V.~Savina, \emph{A reflection formula for the Helmholtz equation
with the Neumann condition}, Comput. Math. Math. Phys. \textbf{39}  (1999),
no. 4, 652-660.
\bibitem{fundsol}
T.V.Savina,  \emph{On splitting up singularities of fundamental
solutions to elliptic equations in $\mathbb{C}^2$}, Cent. Eur.
J. Math., \textbf{5}    (2007), no. 4, 733-740.
\bibitem{bihar}
T.V.~Savina, \emph{On the dependence of the reflection operator on boundary conditions for
biharmonic functions}, J. Math. Anal.  Appl.,  \textbf{370} (2010), 716-725..
\bibitem{shapiro}
H.S.~Shapiro,  \emph{The Schwarz function and its generalization to higher dimensions}, John Wiley and Sons, Inc., 1992.

\bibitem{vekua}
I.N.~Vekua, \emph{New methods for solving elliptic equations}, North
Holland, 1967.





\end{thebibliography}
\end{document}